\documentclass[12pt]{amsart}

\usepackage{amssymb}
\usepackage{amsfonts}
\usepackage{amsmath}
\usepackage{amsthm}
\usepackage[pdftex]{hyperref}
\usepackage{url}
\usepackage{color}
\usepackage{graphicx}

\newcommand{\bq}{\begin{quote}}
\newcommand{\eq}{\end{quote}}
\newcommand{\bi}{\begin{itemize}}
\newcommand{\ei}{\end{itemize}}
\newcommand{\bd}{\begin{description}}
\newcommand{\ed}{\end{description}}
\newcommand{\ben}{\begin{enumerate}}
\newcommand{\een}{\end{enumerate}}
\newcommand{\bbm}{\begin{bmatrix}}
\newcommand{\ebm}{\end{bmatrix}}
\newcommand{\bea}{\begin{eqnarray*}}
\newcommand{\eea}{\end{eqnarray*}}

\def\bF{\mathbb{F}}

\def\KK{\mathsf{K}}

\def\RR{\mathsf{R}}
\def\SS{\mathsf{S}}

\def\XX{\mathsf{X}}
\def\YY{\mathsf{Y}}

\newcommand{\bx}[1]{\mathsf{#1}}     

\def\2G2{\ensuremath{^2{\rm G}_2}}
\def\sl{{\rm{SL}}}

\def\psl{{\rm{PSL}}}
\def\pgl{{\rm{PGL}}}

\def\so{{\rm{SO}}}

\def\a{\alpha}

\def\AXX{\bar{\mathsf{X}}}
\def\TXX{\tilde{\mathsf{X}}}
\def\la{\langle}
\def\ra{\rangle}

\definecolor{darkgreen}{rgb}{0,0.6,0}

\newcommand{\encr}{\ensuremath{\vDash}}

\AtBeginDocument{%
  \def\MR#1{}
  }

\begin{document}

\title{An implementation of the morphisms $\sl_2(\bF) \longrightarrow \sl_2(\KK)  \longrightarrow \XX$}

\author{Alexandre Borovik}
\author{\c{S}\"{u}kr\"{u} Yal\c{c}\i nkaya}


\begin{abstract}
We briefly explain how to implement the morphisms in our paper \cite{BY2025} and provide some examples.
\end{abstract}

\maketitle

\section{Introduction}

In this note, we discuss how to implement  our algorithm in GAP
constructing morphisms
\[
\sl_2(\bF) \rightarrow \sl_2(\KK) \rightarrow \XX.
\]
Here,
\begin{itemize}
\item $\sl_2(\bF)$ is the group of $2\times 2$ matrices of determinant 1 over the field $\bF$ where $\bF$ is a prime field of odd characteristic; 
\item $\XX$ is a black box group encrypting $\sl_2(\bF)$;
\item $\sl_2(\KK)$ is the group of $2\times 2$ matrices of determinant 1 over a black box field $\KK$ encrypting $\bF$, which is constructed inside $\XX$ as presented in \cite{borovik18.540}.
\end{itemize}

Our GAP code is available on
\begin{center}
https://github.com/sukru-yalcinkaya/SL2Morphisms.
\end{center}

Our GAP code is designed to work over prime fields. Our primary focus is to implement our algorithm for the black box groups encrypting $\sl_2$ defined over very big fields so, for simplicity in coding, we assume that the size of the underlying field is at least 13.

\section{About the construction of the group $\pgl_2$}

Let $\XX \encr \sl_2(\bF)$ where $\bF$ is an unknown finite field of odd characteristic. We assume that a set of generators of $\XX$ is given. We work inside the group $\AXX \encr \psl_2(\bF)$ by using the following equivalence of strings in $\XX$:
\begin{equation}\label{eq:psl2}
\bx{x} \equiv \bx{y} \iff \bx{x}\cdot \bx{y}^{-1} \in Z(\XX).
\end{equation}
First, we construct two tori $\SS$ and $\RR$ in $\XX$ where a diagonal automorphism $\bx{\a}$ of $\AXX$ centralizes $\SS$ and inverts $\RR$. To construct a random element in $\AXX$, we consider the elements of the following form:
\[
\bx{x}=s_1r_1\cdots s_kr_k
\]
where $k$ is a reasonably sized random natural number and $s_i$'s and $r_i$'s are random elements from $\SS$ and $\RR$, respectively. We can compute the image of the diagonal automorphism $\a$ for $\bx{x}$:
\[
\bx{x}^\a = s_1r_1^{-1}\cdots s_kr_k^{-1}.
\]

To construct group $\pgl_2(\bF)$, we first consider the diagonal embedding:
\[
\TXX=(\XX, \XX)=\{(\bx{x}, \bx{x}^\a) \mid \bx{x}\in \XX\}.
\]
Clearly, the diagonal automorphism $\a$ interchanges the components of $\TXX$, that is, $(\bx{x}, \bx{x}^\a)^\a=(\bx{x}^\a, \bx{x})$. Finally, $\YY=\TXX \rtimes \la \a \ra \encr \pgl_2(\bF)$ where Equation \ref{eq:psl2} is used for checking whether a group element represents the identity element or not.

In our GAP code, to construct the black box group encrypting $\pgl_2$ by using the black box group $\XX$, we use the following function for the setup:
\[
{\sf{SetUpForPGL2}}(``S",``Eo")
\]
where $S$ is a generating set for $\XX$ and $Eo$ is the odd part of its exponent. This function outputs
\bi
\item The list ``S".
\item A list of elements from a centralizer of an involution, say $i$, inverted by a diagonal automorphism.
\item A list of semisimple elements generating a torus centralized by the same diagonal automorphism.
\item The involution $i$.
\ei

We consider the elements of $\YY$ as follows: 
\[
(\bx{x}, \bx{x}^\a, 0) \mbox{ or } (\bx{x}, \bx{x}^\a, 1)
\]
where the elements of the form $(\bx{x}, \bx{x}^\a, 0)$ belong to the coset $\TXX$ and  the elements of the form $(\bx{x}, \bx{x}^\a, 1)$ belong to the coset $\TXX \a$. The group multiplication in $\YY$ is the usual multiplication in semidirect product of two groups and it is as follows: 
\begin{itemize}
\item $(\bx{x}, \bx{x}^\a, 0) \circ (\bx{y}, \bx{y}^\a, 0) = (\bx{x}\bx{y}, \bx{x}^\a \bx{y}^\a,0).$
\item $(\bx{x}, \bx{x}^\a, 0) \circ  (\bx{y}, \bx{y}^\a, 1) = (\bx{x}\bx{y}, \bx{x}^\a \bx{y}^\a,1).$
\item $(\bx{x}, \bx{x}^\a, 1) \circ  (\bx{y}, \bx{y}^\a, 0) = (\bx{x}\bx{y}^\a, \bx{x}^\a \bx{y},1).$
\item $(\bx{x}, \bx{x}^\a, 1) \circ  (\bx{y}, \bx{y}^\a, 1) = (\bx{x}\bx{y}^\a, \bx{x}^\a \bx{y},0).$
\end{itemize}

\section{Before running the main algorithm}

We need to perform two preprocessing steps before we run our main algorithm.

The first one is {\sf{ToolBoxSL2}}. The function returns all the necessary tools to work within the black box group $\XX$. It can be run as
\[
{\sf{ToolBoxSL2}}(``S", ``E")
\]
where $S$ is a generating set for $\XX$ and $E$ is an any exponent. Its output is the following list.

\bi
\item The output of the function {\sf{SetUpForPGL2}}.
\item A list of elements considered to be a generating set for the semidirect product isomorphic to $\pgl_2$.
\item Three commuting involutions forming the vertices of the projective plane.
\item An element of order 3 permuting the three commuting involutions.
\item A unity element on the corresponding coordinate axes.
\item A generating set for the centralizer of a fixed involution (item 3) which is a vertex determining the projective plane --- the black box field is constructed on the corresponding axis.
\item The projective point (an involution) which serves as 0 in the black box field.
\item The projective point (an involution) corresponding to the homogenous coordinate (1,1,1).
\item Odd part of the exponent of the group generated by the list $S$.
\item Binary representation of the odd part of the exponent.
\item An element of order 4 whose square is the fixed involution.
\item Identity of the group generated by the list $S$.
\ei

Secondly, we construct the change of basis matrix which is used to transform the elements of $\so_3^\sharp$ to the elements of $\so_3^\flat$, see \cite{BY2025} for the definitions of $\so_3^\sharp$ and $\so_3^\flat$.  This function is 
\[
{\sf{SharpVsFlat}}(``TB")
\]
where $``TB"$ is an output of the function {\sf{ToolBoxSL2}}.

\section{How to run our GAP code}

We show over an example how our code should be run in GAP. Consider a black box group $\XX \encr \sl_2(997)$. Since the groups $\sl_2(997)$ exist in GAP Library, we can take its generators and its exponent as follows.

\vskip.3 cm
\begin{center}
\includegraphics[scale=0.16]{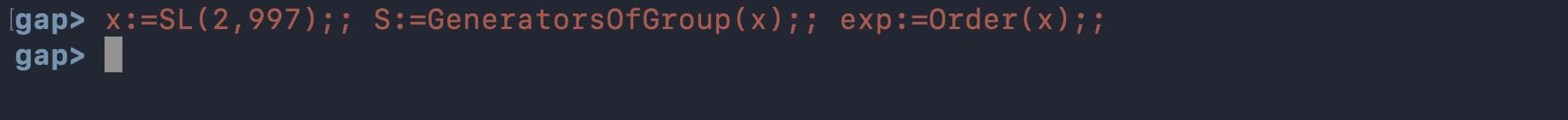}
\end{center}
\vskip.3 cm

We treat the group $\XX$ as a black box group. We construct our tool box for $\XX$ and the change of basis matrix from $\so_3^\sharp$ to $\so_3^\flat$ as follows. The construction of the change of basis matrix involves a computation of a square root of a black box field element and its running time may get lengthy due to some unlucky choices of random elements. Therefore, we provide some information on the screen as the algorithm runs. 

\vskip.3 cm
\begin{center}
\includegraphics[scale=0.16]{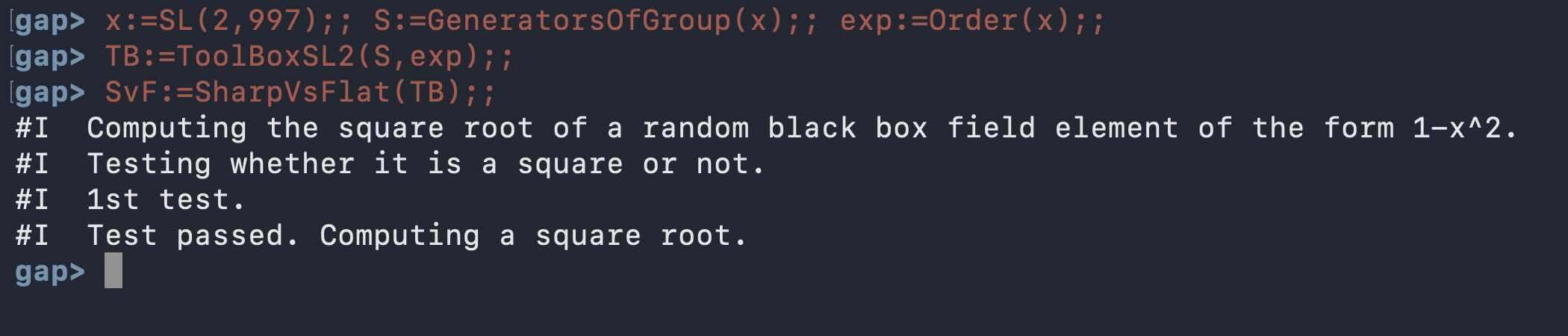}
\end{center}
\vskip.3 cm

Now, we show an example of how an image of random element is found and some simple checks showing that the corresponding images have same orders. 

The input group element is an ordinary element of $\sl_2$ in GAP format. Then, we first represent the $2\times 2$ matrix in GAP as a $2\times 2$ matrix over our black box field by finding the images of the entries of the input matrix in the black box field. Then, we decompose the $2\times 2$ matrix over the black box field $\KK$ as a product of unipotent elements as explained in our paper \cite{BY2025}. After that, we write each unipotent as a product of two involutions and find the images of these involutions by constructing their images in each of the groups below. 
\[
\sl_2(\KK) \rightarrow \so_3^\flat(\KK) \rightarrow \so_3^\sharp(\KK) \rightarrow \XX
\]
As the algorithm runs, we provide the relevant information on the screen. At the end, we take the corresponding component as an output.

Notice that the orders of the corresponding elements are same.
\vskip.3 cm
\begin{center}
\includegraphics[scale=0.16]{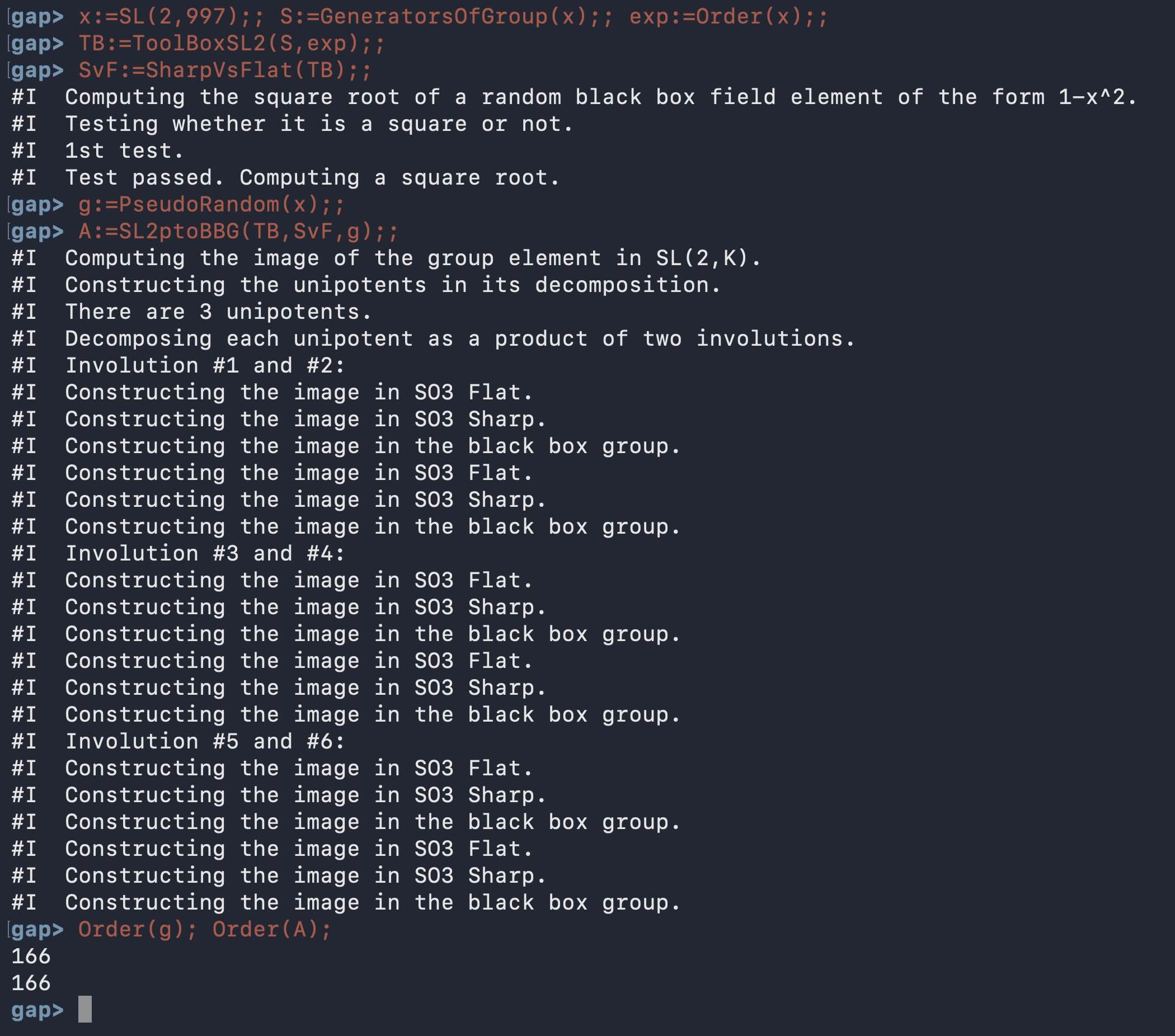}
\end{center}
\vskip.3 cm

For the correctness of the algorithm, it is clear that checking the orders of some random elements and their images is not enough. A more rigorous approach is to examine the Chevalley Commutator Relations for various Chevalley generators. As the constructions of the images of certain Chevalley generators are identical to the example above and the corresponding relations are straightforward to verify, we skip presenting such an example as it would take up too much space in this note. Instead, we present another image of a random element and some checks of the orders of the corresponding elements.

\vskip.3 cm
\begin{center}
\includegraphics[scale=0.16]{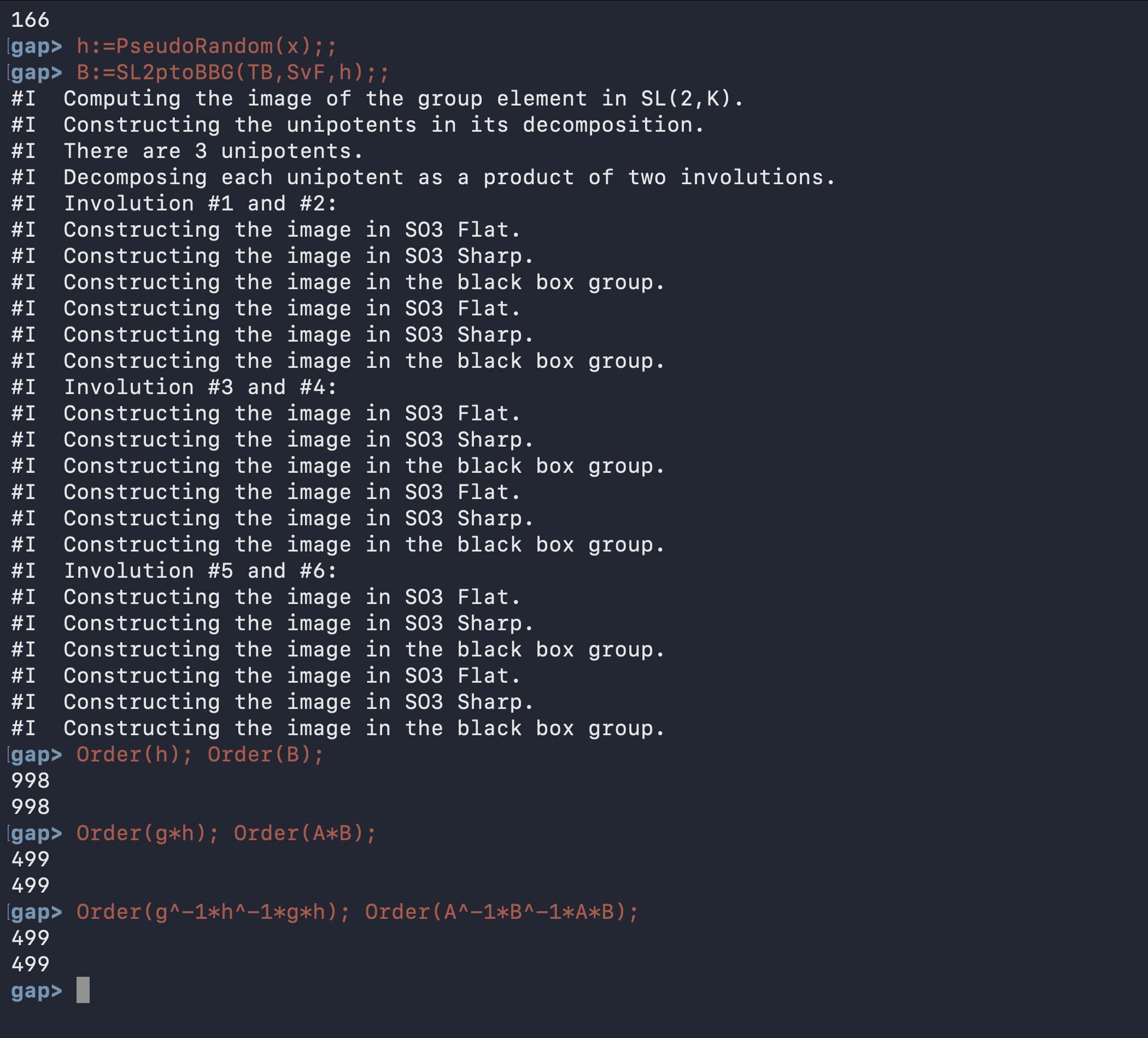}
\end{center}
\vskip.3 cm

In our paper \cite{BY2025}, we presented an algorithm constructing the inverse of this morphism. An implementation of this inverse morphism will be published later in the same GitHub repository.

\end{document}